\documentclass[11pt,a4paper]{article}
\usepackage{amsmath,amssymb,amsthm,mathtools,bm}
\usepackage[margin=1in]{geometry}
\usepackage{microtype}
\usepackage{enumitem}
\usepackage{hyperref}
\usepackage{mathrsfs}
\usepackage{authblk}
\hypersetup{colorlinks=true,linkcolor=blue,citecolor=blue,urlcolor=blue}

\newtheorem{theorem}{Theorem}[section]
\newtheorem{proposition}[theorem]{Proposition}
\newtheorem{lemma}[theorem]{Lemma}
\newtheorem{corollary}[theorem]{Corollary}
\theoremstyle{definition}
\newtheorem{definition}[theorem]{Definition}
\newtheorem{remark}[theorem]{Remark}
\newtheorem{assumption}[theorem]{Assumption}
\newtheorem{problem}[theorem]{Problem}

\newcommand{\R}{\mathbb R}
\newcommand{\ip}[2]{\langle #1,#2\rangle}

\newcommand{\BV}{\mathrm{BV}}

\title{\textbf{A Double-Step Rothe Scheme for Navier--Stokes Inclusions}}
\author{Paweł Szafraniec}
\affil{Faculty of Mathematics and Computer Science \\
	Jagiellonian University in Krak\'ow \\
	ul. \L{}ojasiewicza 6, 30-348 Krak\'ow, Poland}

\begin{document}
\maketitle

\begin{abstract}
In this paper we deal with a double-step Rothe approximation of an evolutionary incompressible Navier--Stokes inclusion with a multivalued term generated by the Clarke subdifferential of a locally Lipschitz functional. The problem is considered in dimensions $d=2$ and $d=3$. We prove solvability of the discrete problems, derive dimension-dependent a priori estimates, and show that the corresponding Rothe functions converge to a weak solution.  Moreover, under additional regularity and relaxed monotonicity, we obtain temporal error estimates and BDF extrapolation gives a second-order estimate.
\end{abstract}



\section{Introduction}

We study an evolutionary incompressible Navier--Stokes inclusion of the form
\[
 u'(t)+Au(t)+B(u(t),u(t))+\gamma^*\xi(t)=f(t),
 \qquad
 \xi(t)\in\partial J(\gamma \gamma(t)),
\]
where $\partial J$ denotes the Clarke subdifferential of a locally Lipschitz functional. The time derivative is approximated by the two-step backward differentiation formula, preceded by one backward Euler step. Our purpose is to establish convergence of the corresponding Rothe approximations in dimensions two and three and, under additional assumptions, to derive rates of convergence.

Rothe methods for evolution inclusions with Clarke subdifferential terms have been considered in a number of papers. In particular, the double-step scheme for parabolic Clarke inclusions was studied in \cite{BartoszSzafraniecZhao2021}, while a second-order evolutionary problem was treated in \cite{BartoszSzafraniec2024}. A one-step Rothe approximation for evolutionary Navier--Stokes hemivariational problems was analyzed in \cite{FangHanMigorskiSofonea2016}. For the classical incompressible Navier--Stokes equations, stability and convergence of the two-step BDF method were investigated by Emmrich \cite{Emmrich2004}. The present paper combines these two ingredients, namely the Navier--Stokes convection and the nonsmooth multivalued term, in a double-step semidiscrete scheme.

The second of the paper is devoted to error estimates. Related estimates for parabolic hemivariational inclusions can be found in \cite{Kalita2012}, and second-order BDF estimates for the incompressible Navier--Stokes equations in \cite{EmmrichError2004}.For the force approximation used in the convergence analysis, we prove an $O(\tau)$ upper error bound. Replacing the force approximation by 
\[
 \widehat f^n=\frac32 f^n-\frac12 f^{n-1}
\]
removes the first-order consistency defect and gives a $O(\tau^2)$ estimate whenever the starting value is smooth enough.

The paper is organized as follows. Section~2 contains notation, assumptions and the formulation of the continuous problem. In Section~3 we introduce the double-step Rothe problem and prove its solvability together with the required a priori estimates. Section~4 is devoted to compactness and convergence. The error estimates are given in Section~5. 

\section{Preliminaries and problem formulation}

Let $\Omega\subset\R^d$, $d\in\{2,3\}$, be a bounded domain with Lipschitz boundary. We introduce the spaces
\[
V=\{v\in H_0^1(\Omega)^d:\operatorname{div}v=0\},
\qquad
H=\overline{V}^{L^2(\Omega)^d}.
\]
As usual, $H$ is identified with its dual. Hence
\[
V\hookrightarrow H\hookrightarrow V^*,
\]
and the embedding $V\hookrightarrow H$ is compact. In what follows, $(\cdot,\cdot)$ and $|\cdot|$ stand for the scalar product and norm in $H$, while $\|\cdot\|$ denotes the norm in $V$. If $i:V\to H$ is the embedding, we put $I=i^*i:V\to V^*$. Thus
\[
\ip{Iv}{w}=(v,w)_H\qquad(v,w\in V).
\]

Let
\[
b(u,v,w)=\int_\Omega (u\cdot\nabla)v\cdot w\,dx,
\qquad
\ip{B(u,v)}{w}=b(u,v,w).
\]
For $u,v,w\in V$,
\[
b(u,v,w)=-b(u,w,v),\qquad b(u,v,v)=0.
\]
\noindent
For $d=2$, the Ladyzhenskaya inequality gives
\begin{equation}\label{eq:2dconv}
\|B(u,u)\|_{V^*}\le C|u|\,\|u\|.
\end{equation}
For $d=3$, interpolation between $L^2(\Omega)^3$ and $L^6(\Omega)^3$ gives
\begin{equation}\label{eq:3dL4}
\|u\|_{L^4(\Omega)}\le C|u|^{1/4}\|u\|^{3/4},
\end{equation}
and hence
\begin{equation}\label{eq:3dconv}
\|B(u,u)\|_{V^*}\le C|u|^{1/2}\|u\|^{3/2}.
\end{equation}
Consequently, for  $u\in L^\infty(0,T;H)\cap L^2(0,T;V)$ we have
\[
B(u,u)\in L^2(0,T;V^*)\quad(d=2),
\]
and
\[
B(u,u)\in L^{4/3}(0,T;V^*)\quad(d=3).
\]

Let $\nu>0$ and define $A:V\to V^*$ by
\[
\ip{Au}{v}=\nu\int_\Omega \nabla u:\nabla v\,dx.
\]
Then
\begin{equation}\label{eq:Aest}
\ip{Au}{u}=\nu\|u\|^2,
\qquad
\|Au\|_{V^*}\le\nu\|u\|.
\end{equation}

Let $U$ be a reflexive Banach space. We assume that $\gamma:V\to U$ is linear and compact and that $J:U\to\R$ is locally Lipschitz. Its Clarke subdifferential, see \cite{Clarke1990} is denoted by $\partial J:U\to2^{U^*}$. We impose the growth condition
\begin{equation}\label{eq:growth}
\|\xi\|_{U^*}\le c_J(1+\|z\|_U)
\qquad\text{for all }z\in U,\ \xi\in\partial J(z).
\end{equation}
When uniqueness in the two-dimensional case or the error estimates of Section~\ref{sec:error} are considered, we additionally assume the relaxed monotonicity condition
\begin{equation}\label{eq:relmono}
\ip{\xi_1-\xi_2}{z_1-z_2}
\ge -m_J\|z_1-z_2\|_U^2
\end{equation}
for $z_i\in U$ and $\xi_i\in\partial J(z_i)$.

We present an Ehrling-type estimate which will be used several times.

\begin{lemma}[Ehrling]\label{lem:ehrling}
Under notation above, for every $\varepsilon>0$ there exists $C_\varepsilon>0$ such that
\begin{equation}\label{eq:ehrling}
\|\gamma v\|_U\le \varepsilon\|v\|+C_\varepsilon|v|
\qquad\text{for every }v\in V.
\end{equation}
\end{lemma}

Combining Lemma~\ref{lem:ehrling} with the growth assumption on $\partial J$, we obtain the following estimate: for every $\eta>0$ there is $C_\eta>0$ such that
\begin{equation}\label{eq:Jenergy}
|\ip{\xi}{\gamma v}|
\le \eta\|v\|^2+C_\eta(1+|v|^2)
\end{equation}
for all $v\in V$ and $\xi\in\partial J(\gamma v)$. In fact, \eqref{eq:growth} yields
\[
|\ip{\xi}{\gamma v}|
\le c_J(1+\|\gamma v\|_U)\|\gamma v\|_U,
\]
and then the result follows from \eqref{eq:ehrling} and Young's inequality.

\begin{definition}\label{defininitionMpq}
	
	Let $1\le q<\infty$ and let $X$ be a Banach space. For a $w:[0,T]\to X$ put
	\[
	|w|_{\BV^q(0,T;X)}^q
	=\sup_{\pi}\sum_{j=1}^m\|w(t_j)-w(t_{j-1})\|_X^q,
	\]
	where the supremum is taken over all finite partitions $0=t_0<\cdots<t_m=T$. For Banach spaces $X_1\hookrightarrow X_3$ we set
	\[
	M^{p,q}(0,T;X_1,X_3)
	=L^p(0,T;X_1)\cap\BV^q(0,T;X_3),
	\]
	with
	\[
	\|w\|_{M^{p,q}}
	=\|w\|_{L^p(0,T;X_1)}+|w|_{\BV^q(0,T;X_3)}.
	\]
\end{definition}

The next result is crucial for the convergence estimates, see \cite{Kalita2012}.

\begin{lemma}\label{lemmaKalita}
	 If $X_1$ is reflexive, $X_1\hookrightarrow X_2$ compactly, and $X_2\hookrightarrow X_3$ continuously, then bounded subsets of $M^{2,2}(0,T;X_1,X_3)$ are relatively compact in $L^2(0,T;X_2)$.
\end{lemma}
\

We ready to formulate the classical problem.

\begin{problem}[Problem classical]
Let $f\in L^2(0,T;V^*)$ and $u_0\in H$. We will call 
\begin{equation}\label{eq:problem}
\left\{
\begin{aligned}
u'(t)+Au(t)+B(u(t),u(t))+\gamma^*\xi(t)&=f(t),\\
\xi(t)&\in\partial J(\gamma u(t)),\\
u(0)&=u_0.
\end{aligned}
\right.
\end{equation}
	zawartość...
\end{problem}

For convenience, put
\begin{equation}\label{eq:qd}
q_d=\begin{cases}2,&d=2,\\4/3,&d=3.\end{cases}
\end{equation}

\begin{definition}\label{def:weak}
A weak solution of \eqref{eq:problem} is a pair $(u,\xi)$ such that
\[
u\in L^\infty(0,T;H)\cap L^2(0,T;V),
\qquad
u'\in L^{q_d}(0,T;V^*),
\]
\[
\xi\in L^2(0,T;U^*),
\qquad
\xi(t)\in\partial J(\gamma u(t))\quad\text{for a.e. }t,
\]
and the first relation in \eqref{eq:problem} is satisfied in $V^*$ for a.e. $t\in(0,T)$.
\end{definition}

\section{The double-step Rothe scheme}

Let $N\in\mathbb N$, $\tau=T/N$, and $t_n=n\tau$ for $n=0,\ldots,N$. For the right-hand side we put
\[
f^n=\frac1\tau\int_{t_{n-1}}^{t_n}f(t)\,dt.
\]
We choose $u^0_\tau\in V$ satisfying
\begin{equation}\label{eq:initapprox}
u^0_\tau\to u_0\quad\text{strongly in }H.
\end{equation}
Such a sequence exists by the density of $V$ in $H$.

The first discrete problem is given by the backward Euler formula
\begin{equation}\label{eq:BE}
\frac{u^1-u^0}{\tau}+Au^1+B(u^1,u^1)+\gamma^*\xi^1=f^1,
\qquad \xi^1\in\partial J(\gamma u^1).
\end{equation}
For $n=2,\ldots,N$ we set
\[
D_2u^n=\frac1\tau\left(\frac32u^n-2u^{n-1}+\frac12u^{n-2}\right),
\]
and consider the double-step problem
\begin{equation}\label{eq:BDF}
D_2u^n+Au^n+B(u^n,u^n)+\gamma^*\xi^n=f^n,
\qquad \xi^n\in\partial J(\gamma u^n).
\end{equation}

\begin{lemma}\label{lem:Bcompact}
For $d=2,3$, the map $v\mapsto B(v,v)$ is completely continuous from $V$ endowed with its weak topology into $V^*$ endowed with its norm topology. Consequently $A+B(\cdot,\cdot)$ is pseudomonotone.
\end{lemma}

\begin{proof}
Let $v_k\rightharpoonup v$ in $V$. Since $H_0^1(\Omega)^d\hookrightarrow L^4(\Omega)^d$ is compact for $d\in\{2,3\}$, it follows that $v_k\to v$ in $L^4(\Omega)^d$. For $w\in V$, $\|w\|\le1$, the skew-symmetry of $b$ gives
\[
\begin{aligned}
|b(v_k,v_k,w)-b(v,v,w)|
&\le |b(v_k-v,v_k,w)|+|b(v,v_k-v,w)|\\
&= |b(v_k-v,w,v_k)|+|b(v,w,v_k-v)|\\
&\le C\|v_k-v\|_{L^4}
\bigl(\|v_k\|_{L^4}+\|v\|_{L^4}\bigr)\|w\|.
\end{aligned}
\]
The last expression tends to zero. Hence $B(v_k,v_k)\to B(v,v)$ in $V^*$. Thus $B(\cdot,\cdot)$ is completely continuous on the diagonal. Since $A$ is bounded, linear, monotone and hemicontinuous, the operator $A+B(\cdot,\cdot)$ is pseudomonotone.
\end{proof}

Fisrt we prove existence of discrete problem at each step

\begin{theorem}[Existence of the dicrete problem]\label{thm:discreteexist}
There exists $\tau_{\rm sol}>0$, depending only on $\nu$, $c_J$, $\gamma$ and the constants in the Ehrling estimate, such that for every $0<\tau\le\tau_{\rm sol}$ the problems \eqref{eq:BE} and \eqref{eq:BDF} admit solutions in both dimensions $d=2$ and $d=3$.
\end{theorem}

\begin{proof}
For $n=2,\ldots,N$ we set $\mathcal T_n:V\to2^{V^*}$ by
\[
\mathcal T_n(v)=\frac32v+\tau Av+\tau B(v,v)+\tau\gamma^*\partial J(\gamma v).
\]
Then \eqref{eq:BDF} is equivalent to
\[
2u^{n-1}-\frac12u^{n-2}+\tau f^n\in\mathcal T_n(v).
\]
We apply the surjectivity theorem for coercive pseudomonotone multifunctions, cf. \cite[Theorem~3.70]{MigorskiOchalSofonea2013}. For every $v\in V$, the set $\partial J(\gamma v)$ is nonempty, convex and weakly compact in $U^*$. Hence $\gamma^*\partial J(\gamma v)$ and $\mathcal T_n(v)$ have nonempty, convex and weakly compact values in $V^*$. By \eqref{eq:growth}, \eqref{eq:Aest} and
\[
\|B(v,v)\|_{V^*}\le C\|v\|^2\qquad(d\in\{2,3\}),
\]
the multifunction $\mathcal T_n$ is bounded on bounded subsets of $V$. The Clarke subdifferential is strong-to-weak upper semicontinuous on bounded sets. Since $\gamma$ is compact, the multifunction $v\mapsto\gamma^*\partial J(\gamma v)$ is pseudomonotone. Together with Lemma~\ref{lem:Bcompact} and the sum theorem, this proves that $\mathcal T_n$ is pseudomonotone.

Let $v\in V$ and $w\in\mathcal T_n(v)$. Then, there exists $\eta\in\partial J(\gamma v)$ such that
\[
\ip{w}{v}=\frac32|v|^2+\tau\nu\|v\|^2+\tau\ip{\eta}{\gamma v},
\]
since $b(v,v,v)=0$. Applying \eqref{eq:Jenergy} with the coefficient of $\|v\|^2$ chosen as $\nu/2$:
\[
\ip{w}{v}
\ge \left(\frac32-C_E\tau\right)|v|^2
+\frac{\tau\nu}{2}\|v\|^2-C_E\tau,
\]
where $C_E>0$ depends only on $c_J$, $\gamma$ and the constant in the Ehrling inequality. We choose $\tau_{\rm sol}>0$ such that $C_E\tau_{\rm sol}\le3/4$. We obtain
\[
\ip{w}{v}\ge \frac34|v|^2+\frac{\tau\nu}{2}\|v\|^2-C_E\tau,
\]
and thus $\mathcal T_n$ is coercive for every $0<\tau\le\tau_{\rm sol}$ and the classical surjectivity theorem gives a solution of \eqref{eq:BDF}. The proof for \eqref{eq:BE} is the same, with $I$ in place of $\frac32 I$. If necessary, $\tau_{\rm sol}$ is decreased so that the coefficient of the $H$-term remains positive.
\end{proof}

We shall repeatedly use the identity
\begin{equation}\label{eq:Gidentity}
\left(\frac32a-2b+\frac12c,a\right)
=\frac14\Bigl(|a|^2+|2a-b|^2-|b|^2-|2b-c|^2+|a-2b+c|^2\Bigr)
\end{equation}
which is the classical $G$-stability identity for BDF; see \cite{Dahlquist1978,NevanlinnaOdeh1981} and, for the incompressible Navier--Stokes equations, \cite{Emmrich2004}.

\begin{proposition}\label{prop:energy}
There exists $\tau_{\rm en}>0$, depending on $\nu$, $c_J$, $\gamma$ and the constants in the Ehrling estimate, such that for $0<\tau\le\tau_{\rm en}$ there is a constant $C_T$ independent of $\tau$ such that
\begin{equation}\label{eq:energy}
\max_{0\le n\le N}|u^n|^2
+\sum_{n=2}^N|u^n-2u^{n-1}+u^{n-2}|^2
+\tau\sum_{n=1}^N\|u^n\|^2\le C_T.
\end{equation}
Moreover,
\begin{equation}\label{eq:xiest}
\tau\sum_{n=1}^N\|\xi^n\|_{U^*}^2\le C_T.
\end{equation}
\end{proposition}

\begin{proof}
We test \eqref{eq:BE} with $u^1$. Since $b(u^1,u^1,u^1)=0$, we obtain
\[
\frac12\bigl(|u^1|^2-|u^0|^2+|u^1-u^0|^2\bigr)
+\tau\nu\|u^1\|^2
+\tau\ip{\xi^1}{\gamma u^1}
=\tau\ip{f^1}{u^1}.
\]
Fix $0<\eta<\nu/4$. By \eqref{eq:Jenergy} and Young's inequality,
\[
|\ip{\xi^1}{\gamma u^1}|
\le \eta\|u^1\|^2+C_\eta(1+|u^1|^2),
\]
\[
|\ip{f^1}{u^1}|
\le \eta\|u^1\|^2+C_\eta\|f^1\|_{V^*}^2.
\]
Hence
\begin{equation}\label{eq:firstenergy}
|u^1|^2-|u^0|^2+|u^1-u^0|^2+c\tau\|u^1\|^2
\le C\tau\bigl(1+|u^1|^2+\|f^1\|_{V^*}^2\bigr).
\end{equation}

For $n\ge2$, we test \eqref{eq:BDF} with $u^n$ and multiply by $\tau$. By \eqref{eq:Gidentity}, we get
\[
\begin{aligned}
&\frac14\Bigl(|u^n|^2+|2u^n-u^{n-1}|^2-|u^{n-1}|^2-|2u^{n-1}-u^{n-2}|^2\\
&\hspace{27mm}+|u^n-2u^{n-1}+u^{n-2}|^2\Bigr)+\tau\nu\|u^n\|^2\\
&\qquad +\tau\ip{\xi^n}{\gamma u^n}=\tau\ip{f^n}{u^n}.
\end{aligned}
\]
Arguing as for the first step, we get
\[
\begin{aligned}
&|u^n|^2+|2u^n-u^{n-1}|^2-|u^{n-1}|^2-|2u^{n-1}-u^{n-2}|^2+|u^n-2u^{n-1}+u^{n-2}|^2+c\tau\|u^n\|^2\\
&\qquad\le C\tau\bigl(1+|u^n|^2+\|f^n\|_{V^*}^2\bigr).
\end{aligned}
\]

Since $|u^n|^2\le |u^n|^2+|2u^n-u^{n-1}|^2.$, we sum the above inequalities and apply the discrete Gronwall lemma. Taking $\tau_{\rm en}>0$ sufficiently small to absorb the term $C\tau E_n$, we obtain \eqref{eq:energy}. Moreover, there are constants $C_0,C_1>0$, independent of $T$ and $\tau$, such that
\[
C_T\le C_0 e^{C_1T}
\left(1+|u^0|^2+\|f\|_{L^2(0,T;V^*)}^2\right).
\]
The constants $C_0$ and $C_1$ depend on $\nu$, $c_J$, $\gamma$ and the Ehrling constant. 

To prove \eqref{eq:xiest}, we use \eqref{eq:growth} and \eqref{eq:ehrling}. For a fixed $\varepsilon>0$,
\[
\|\xi^n\|_{U^*}^2
\le C\bigl(1+\|\gamma u^n\|_U^2\bigr)
\le C\bigl(1+\|u^n\|^2+|u^n|^2\bigr).
\]
After multiplication by $\tau$ and summation, \eqref{eq:xiest} follows from \eqref{eq:energy}.
\end{proof}

\begin{remark}\label{rem:taustar}
In what follows we put 
\[
\tau_*:=\min\{\tau_{\rm sol},\tau_{\rm en}\}.
\]
All subsequent discrete estimates are understood for $0<\tau\le\tau_*$. The number $\tau_*$ depends only on the structural constants, whereas $C_T$ may depend on $T$ through the Gronwall inequaltiy
\end{remark}

We introduce the piecewise constant function
\[
d_\tau u(t)=
\begin{cases}
(u^1-u^0)/\tau,&t\in(0,\tau],\\
D_2u^n,&t\in(t_{n-1},t_n],\quad n\ge2.
\end{cases}
\]
We also set $\delta_n=u^n-u^{n-1}$.

\begin{proposition}\label{prop:d2}
If $d=2$, then
\begin{equation}\label{eq:d2der}
\|d_\tau u\|_{L^2(0,T;V^*)}\le C,
\end{equation}
and
\begin{equation}\label{eq:d2inc}
\sum_{n=1}^N\|\delta_n\|_{V^*}^2\le C\tau.
\end{equation}
\end{proposition}

\begin{proof}
From the discrete equations and \eqref{eq:2dconv},
\[
\|D_2u^n\|_{V^*}
\le \|f^n\|_{V^*}+\nu\|u^n\|+C|u^n|\|u^n\|+\|\gamma\|\|\xi^n\|_{U^*}
\]
for $n\ge2$, and the analogous estimate holds for $(u^1-u^0)/\tau$. Hence \eqref{eq:d2der} follows from Proposition~\ref{prop:energy} and \eqref{eq:xiest}.

For the first step, \eqref{eq:BE} gives
\[
\|\delta_1\|_{V^*}
\le\tau\Bigl(\|f^1\|_{V^*}+\nu\|u^1\|+C|u^1|\|u^1\|+\|\gamma\|\|\xi^1\|_{U^*}\Bigr).
\]
Since
\[
\tau\|f^1\|_{V^*}^2+\tau\|u^1\|^2+\tau\|\xi^1\|_{U^*}^2\le C
\]
and $|u^1|\le C$, it follows that
\begin{equation}\label{eq:delta1d2}
\|\delta_1\|_{V^*}\le C\tau^{1/2},
\qquad
\|\delta_1\|_{V^*}^2\le C\tau.
\end{equation}
For $n\ge2$,
\begin{equation}\label{eq:deltaRec}
\delta_n=\frac23\tau D_2u^n+\frac13\delta_{n-1}.
\end{equation}
Hence
\[
\|\delta_n\|_{V^*}^2
\le\frac89\tau^2\|D_2u^n\|_{V^*}^2+\frac29\|\delta_{n-1}\|_{V^*}^2.
\]
Summation of the recurrence, together with \eqref{eq:delta1d2} and \eqref{eq:d2der}, gives
\[
\sum_{n=1}^N\|\delta_n\|_{V^*}^2
\le C\tau+C\tau^2\sum_{n=2}^N\|D_2u^n\|_{V^*}^2
\le C\tau.
\]
\end{proof}

\begin{proposition}\label{prop:d3}
If $d=3$, then
\begin{equation}\label{eq:d3der}
\tau\left\|\frac{u^1-u^0}{\tau}\right\|_{V^*}^{4/3}
+\tau\sum_{n=2}^N\|D_2u^n\|_{V^*}^{4/3}\le C,
\end{equation}
and
\begin{equation}\label{eq:d3inc}
\sum_{n=1}^N\|\delta_n\|_{V^*}^{4/3}\le C\tau^{1/3}.
\end{equation}
\end{proposition}

\begin{proof}
For $n\ge2$, \eqref{eq:3dconv} gives
\[
\|D_2u^n\|_{V^*}
\le \|f^n\|_{V^*}+\nu\|u^n\|
+C|u^n|^{1/2}\|u^n\|^{3/2}
+\|\gamma\|\|\xi^n\|_{U^*}.
\]
We raise the inequality to the power $4/3$, multiply by $\tau$, and sum over $n$. This gives
\[
\begin{aligned}
\tau\sum_{n=2}^N\|D_2u^n\|_{V^*}^{4/3}
&\le C\tau\sum_{n=2}^N\|f^n\|_{V^*}^{4/3}
+C\tau\sum_{n=2}^N\|u^n\|^{4/3}\\
&\quad+C\tau\sum_{n=2}^N|u^n|^{2/3}\|u^n\|^2
+C\tau\sum_{n=2}^N\|\xi^n\|_{U^*}^{4/3}\le C.
\end{aligned}
\]
The first, second and fourth terms are bounded straightforwardly by Proposition~\ref{prop:energy}. Moreover, to deal with the third term, we have again using \ref{prop:energy}
\[
\tau\sum_{n=2}^N|u^n|^{2/3}\|u^n\|^2
\le \Bigl(\max_n|u^n|^2\Bigr)^{1/3}
\tau\sum_{n=2}^N\|u^n\|^2\le C.
\]

For the starting point, \eqref{eq:BE} and \eqref{eq:3dconv} give
\[
\begin{aligned}
\|\delta_1\|_{V^*}
&\le \tau\|f^1\|_{V^*}+\tau\nu\|u^1\|
+C\tau |u^1|^{1/2}\|u^1\|^{3/2}
+\tau\|\gamma\|\|\xi^1\|_{U^*}.
\end{aligned}
\]
The first, second and fourth terms are of order $\tau^{1/2}$. For the convective term we only have the order $\tau^{1/4}$, since $|u^1|\le C$ and $\tau\|u^1\|^2\le C$. Therefore
\begin{equation}\label{eq:delta1d3}
\|\delta_1\|_{V^*}\le C\tau^{1/4},
\qquad
\|\delta_1\|_{V^*}^{4/3}\le C\tau^{1/3}.
\end{equation}
Consequently,
\[
\tau\left\|\frac{\delta_1}{\tau}\right\|_{V^*}^{4/3}
=\tau^{-1/3}\|\delta_1\|_{V^*}^{4/3}\le C,
\]
and \eqref{eq:d3der} follows.

 From \eqref{eq:deltaRec} and the convexity inequality, one can choose $\theta\in(0,1)$ and $C>0$ such that
\[
\|\delta_n\|_{V^*}^{\frac43}
\le C\tau^\frac43\|D_2u^n\|_{V^*}^{\frac43}+\theta\|\delta_{n-1}\|_{V^*}^{\frac43}.
\]
Iterating this inequality and summing with respect to $n$ we obtain
\[
\sum_{n=1}^N\|\delta_n\|_{V^*}^{\frac43}
\le C\|\delta_1\|_{V^*}^{\frac43}
+C\tau^\frac43\sum_{n=2}^N\|D_2u^n\|_{V^*}^{\frac43}.
\]
The first term is of order $\tau^{1/3}$, and the second one is
\[
C\tau^{\frac43-1}\left(\tau\sum_{n=2}^N\|D_2u^n\|_{V^*}^{\frac43}\right)
\le C\tau^{1/3}.
\]
This proves \eqref{eq:d3inc}.
\end{proof}

Let
\[
\bar u_\tau(t)=u^n\quad\text{for }t\in(t_{n-1},t_n],
\qquad \bar u_\tau(0)=u^0.
\]

\begin{proposition}\label{prop:Mq}
For $d\in\{2,3\}$,
\begin{equation}\label{eq:Mqd}
\|\bar u_\tau\|_{M^{2,q_d}(0,T;V,V^*)}\le C.
\end{equation}
Thus the natural compactness class is $M^{2,2}$ in $d=2$ and $M^{2,4/3}$ in $d=3$.
\end{proposition}

\begin{proof}
The $L^2(0,T;V)$ estimate is contained in \eqref{eq:energy}. Put $q=q_d$ and let $0=t_0<\cdots<t_m=T$ be an arbitrary partition. Every nonzero increment $\bar u_\tau(t_j)-\bar u_\tau(t_{j-1})$ is a sum of consecutive differences $\delta_r$. If $I_j$ contains $m_j$ such differences, then
\[
\left\|\sum_{r\in I_j}\delta_r\right\|_{V^*}^{q}
\le m_j^{q-1}\sum_{r\in I_j}\|\delta_r\|_{V^*}^{q}.
\]
The sets $I_j$ are disjoint and $m_j\le N=T/\tau$. Hence
\[
\sum_j\|\bar u_\tau(t_j)-\bar u_\tau(t_{j-1})\|_{V^*}^{q}
\le (T/\tau)^{q-1}\sum_{r=1}^N\|\delta_r\|_{V^*}^{q}.
\]
For $d=2$, \eqref{eq:d2inc} gives a uniform bound. For $d=3$, \eqref{eq:d3inc} gives
\[
(T/\tau)^{1/3}C\tau^{1/3}\le CT^{1/3}.
\]
Taking the supremum over all partitions yields \eqref{eq:Mqd}.
\end{proof}

The following observation allows us to consider the same space for both dimmensions.

\begin{proposition} \label{prop:M22all}
If $d=3$, then the $\bar u_\tau$  is also bounded in $M^{2,2}(0,T;V,V^*)$.
\end{proposition}

\begin{proof}
Since $H\hookrightarrow V^*$ continuously and \eqref{eq:energy} gives $\|\bar u_\tau\|_{L^\infty(0,T;H)}\le C$, there is $M>0$ with
\[
\sup_{t\in[0,T]}\|\bar u_\tau(t)\|_{V^*}\le M.
\]
Hence every partition increment is bounded by $2M$, and therefore
\[
\begin{aligned}
\sum_j\|\bar u_\tau(t_j)-\bar u_\tau(t_{j-1})\|_{V^*}^2
&\le (2M)^{2/3}
\sum_j\|\bar u_\tau(t_j)-\bar u_\tau(t_{j-1})\|_{V^*}^{4/3}\\
&\le C.
\end{aligned}
\]
The result now follows by taking the supremum over all partitions and using Proposition~\ref{prop:Mq}.
\end{proof}

\section{Convergence}\label{sec:compactness}

\begin{proposition}\label{prop:compact}
There exist a subsequence and $u\in L^\infty(0,T;H)\cap L^2(0,T;V)$ such that
\[
\bar u_\tau\rightharpoonup u\quad\text{in }L^2(0,T;V),
\]
\[
\bar u_\tau\overset{*}{\rightharpoonup}u\quad\text{in }L^\infty(0,T;H),
\]
\begin{equation}\label{eq:strongH}
\bar u_\tau\to u\quad\text{strongly in }L^2(0,T;H),
\end{equation}
and hence
\[
\bar u_\tau(t)\to u(t)\quad\text{in }H\quad\text{for a.e. }t \in (0,T)
\]
Moreover,
\begin{equation}\label{eq:strongell}
\gamma\bar u_\tau\to\gamma u\quad\text{strongly in }L^2(0,T;U).
\end{equation}
\end{proposition}

\begin{proof}
The first two convergences follow from \eqref{eq:energy}. Applying Proposition~\ref{prop:Mq} in dimension two and Proposition~\ref{prop:M22all} in dimension three, and then the above compactness theorem with $p=q=2$ to the triple
\[
V\hookrightarrow H\hookrightarrow V^*
\]
we obtain \eqref{eq:strongH}. Passing to a further subsequence, if necessary, gives the pointwise convergence in $H$ for a.e. $t\in(0,T)$.

It remains to prove the last convergence. Put $w_\tau=\bar u_\tau-u$. Then $(w_\tau)$ is bounded in $L^2(0,T;V)$ and $w_\tau\to0$ in $L^2(0,T;H)$. Lemma~\ref{lem:ehrling} gives, for arbitrary $\varepsilon>0$,
\[
\|\gamma w_\tau(t)\|_U
\le\varepsilon\|w_\tau(t)\|+C_\varepsilon|w_\tau(t)|.
\]
Therefore
\[
\|\gamma w_\tau\|_{L^2(0,T;U)}^2
\le2\varepsilon^2\|w_\tau\|_{L^2(0,T;V)}^2
+2C_\varepsilon^2\|w_\tau\|_{L^2(0,T;H)}^2.
\]
Passing to the upper limit as $\tau\to0$, we get
\[
\limsup_{\tau\to0}\|\gamma w_\tau\|_{L^2(0,T;U)}^2\le C\varepsilon^2.
\]
Since $\varepsilon>0$ is arbitrary, this proves \eqref{eq:strongell}. 
\end{proof}

Since $T<\infty$, the convergence in \eqref{eq:xiweak} also yields the weak $L^1(0,T;U^*)$ convergence required in the Aubin--Cellina theorem.

\begin{theorem}\label{thm:Clarke} 
	Let $\bar\xi_\tau(t)=\xi^n$ for $t\in(t_{n-1},t_n]$. From \eqref{eq:xiest}, after passing to a subsequence,
	\begin{equation}\label{eq:xiweak}
		\bar\xi_\tau\rightharpoonup\xi\quad\text{in }L^2(0,T;U^*).
	\end{equation}
The limit satisfies
\[
\xi(t)\in\partial J(\gamma u(t))\qquad\text{for a.e. }t\in(0,T).
\]
\end{theorem}

\begin{proof}
By \eqref{eq:strongell}, after taking a further subsequence,
\[
\gamma\bar u_\tau(t)\to\gamma u(t)\quad\text{strongly in }U
\]
for almost every $t$. Also
\[
\bar\xi_\tau(t)\in\partial J(\gamma\bar u_\tau(t))
\quad\text{for a.e. }t.
\]
The multifunction $\partial J$ has nonempty, convex and weakly compact values on bounded sets and is upper semicontinuous from the strong topology of $U$ to the weak topology of $U^*$. Since $T<\infty$, the convergence in \eqref{eq:xiweak} also yields the weak $L^1(0,T;U^*)$ convergence. Hence the Aubin--Cellina convergence theorem \cite[Theorem~1.4.1]{AubinCellina1984} gives
$\xi(t)\in\partial J(\gamma u(t))$ for a.e. $t\in(0,T)$.
\end{proof}

\begin{lemma}\label{lem:conv}
We have
\[
B(\bar u_\tau,\bar u_\tau)\rightharpoonup B(u,u)
\quad\text{in }L^{q_d}(0,T;V^*),
\]
where $q_d$ is given by \eqref{eq:qd}.
\end{lemma}

\begin{proof}
The boundedness in $L^{q_d}(0,T;V^*)$ follows directly from \eqref{eq:2dconv}, \eqref{eq:3dconv} and \eqref{eq:energy}.  To identify the weak limit, observe first that
\[
u\in L^\infty(0,T;H)\cap L^2(0,T;V).
\]
In dimension three, interpolation between \eqref{eq:strongH} and the uniform bound in $L^2(0,T;L^6(\Omega)^3)$ provides us for every $2\le r<6$ with
\begin{equation}\label{eq:strongLr}
\bar u_\tau\to u\quad\text{strongly in }L^2(0,T;L^r(\Omega)^3).
\end{equation}
Now take $r=3$, put $w_\tau=\bar u_\tau-u$ and fix $v\in V$. By the skew-symmetry of $b$,
\[
b(w_\tau,\bar u_\tau,v)=-b(w_\tau,v,\bar u_\tau),
\qquad
b(u,w_\tau,v)=-b(u,v,w_\tau).
\]
Therefore
\[
\begin{aligned}
|b(\bar u_\tau,\bar u_\tau,v)-b(u,u,v)|
&\le C\|w_\tau\|_{L^3}
\bigl(\|\bar u_\tau\|_{L^6}+\|u\|_{L^6}\bigr)\|v\|.
\end{aligned}
\]
The right-hand side converges to zero in $L^1(0,T)$ by H\"older's inequality and \eqref{eq:strongLr}. The two-dimensional case is proved in the same way, using the compact embedding of $H_0^1(\Omega)$ into $L^r(\Omega)$ for $1\le r<\infty$.

We have therefore proved convergence for functions of the form $\varphi(t)v$, where $\varphi\in C_0^\infty(0,T)$ and $v\in V$. Their finite linear combinations are dense in $L^{q_d'}(0,T;V)$, $q_d'=q_d/(q_d-1)$. Since $B(\bar u_\tau,\bar u_\tau)$ is uniformly bounded in $L^{q_d}(0,T;V^*)$, a density argument completes the proof.
\end{proof}

We next establish the summation-by-parts formula needed to pass to the limit in the discrete derivative. Let $\varphi\in C^3([0,T])$, $\varphi(T)=0$, and set $\varphi_n=\varphi(t_n)$. Put
\[
S_\tau(v,\varphi)
=(u^1-u^0,v)\varphi_1
+\sum_{n=2}^N\left(\frac32u^n-2u^{n-1}+\frac12u^{n-2},v\right)\varphi_n.
\]
After reaggaging the sums, we can write it as
\begin{equation}\label{eq:SBPexact}
\begin{aligned}
S_\tau(v,\varphi)
&=(u^0,v)\left(-\varphi_1+\frac12\varphi_2\right)
+(u^1,v)\left(\varphi_1-2\varphi_2+\frac12\varphi_3\right)\\
&\quad+\sum_{k=2}^{N-2}(u^k,v)
\left(\frac32\varphi_k-2\varphi_{k+1}+\frac12\varphi_{k+2}\right)\\
&\quad+(u^{N-1},v)\left(\frac32\varphi_{N-1}-2\varphi_N\right)
+\frac32(u^N,v)\varphi_N.
\end{aligned}
\end{equation}
Since $\varphi_N=0$, the last term is equal to zero.

\begin{lemma}\label{lem:SBP}
For every $v\in V$ and every $\varphi\in C^3([0,T])$ with $\varphi(T)=0$,
\begin{equation}\label{eq:SBPlimit}
S_\tau(v,\varphi)
\longrightarrow
-\int_0^T(u(t),v)\varphi'(t)\,dt-(u_0,v)\varphi(0).
\end{equation}
\end{lemma}

\begin{proof}
For $2\le k\le N-2$, Taylor's formula at $t_k$ yields
\begin{equation}\label{eq:taylork}
\frac32\varphi_k-2\varphi_{k+1}+\frac12\varphi_{k+2}
=-\tau\varphi'(t_k)+R_k,
\qquad |R_k|\le C\tau^3\|\varphi^{(3)}\|_\infty.
\end{equation}
By \eqref{eq:energy}, $\max_k|u^k|\le C$, and hence the sum of all remainder terms is $O(N\tau^3)=O(\tau^2)$.

At the left endpoint, Taylor's formula gives
\[
-\varphi_1+\frac12\varphi_2
=-\frac12\varphi(0)+O(\tau^2),
\]
and
\[
\varphi_1-2\varphi_2+\frac12\varphi_3
=-\frac12\varphi(0)-\frac32\tau\varphi'(0)+O(\tau^2).
\]
Furthermore, \eqref{eq:delta1d2} and \eqref{eq:delta1d3} imply, respectively in dimensions two and three,
\[
\|u^1-u^0\|_{V^*}\to0.
\]
Since $u^0_\tau\to u_0$ in $H$, it follows that $(u^1,v)\to(u_0,v)$ for every $v\in V$. Consequently, the two initial terms in \eqref{eq:SBPexact} converge to $-(u_0,v)\varphi(0)$.

At the right endpoint, $\varphi_N=0$ and $|\varphi_{N-1}|\le\tau\|\varphi'\|_\infty$. Hence the term containing $u^{N-1}$ tends to zero.

It remains to treat the interior sum. Define the step function
\[
\psi_\tau(t)=\varphi'(t_k)\quad\text{for }t\in(t_{k-1},t_k],\quad 2\le k\le N-2.
\]
The uniform continuity of $\varphi'$ gives $\psi_\tau\to\varphi'$ uniformly on the union of these intervals. Moreover,
\[
-\tau\sum_{k=2}^{N-2}(u^k,v)\varphi'(t_k)
=-\int_{t_1}^{t_{N-2}}(\bar u_\tau(t),v)\psi_\tau(t)\,dt.
\]
The omitted intervals have total length at most $3\tau$, and their contribution is $O(\tau)$. On the remaining interval, \eqref{eq:strongH} and the uniform convergence of $\psi_\tau$ yield
\[
\int_{t_1}^{t_{N-2}}(\bar u_\tau,v)\psi_\tau\,dt
\longrightarrow
\int_0^T(u,v)\varphi'\,dt.
\]
Combining the above limits with the $O(\tau^2)$ Taylor remainder proves \eqref{eq:SBPlimit}.
\end{proof}

\begin{theorem}\label{thm:main}
Let $d\in\{2,3\}$, $u_0\in H$, $f\in L^2(0,T;V^*)$, and assume \eqref{eq:growth}. Let the initial approximations $u^0_\tau\in V$ satisfy \eqref{eq:initapprox}. Let $\tau\downarrow0$ with $0<\tau\le\tau_*$, and let $(u^n,\xi^n)$ be any choice of solutions of the BDF--Rothe scheme at the successive time steps. Then, there exist a subsequence and a weak solution $(u,\xi)$ in the sense of Definition \ref{def:weak} such that
\[
\bar u_\tau\rightharpoonup u\quad\text{in }L^2(0,T;V),
\]
\[
\bar u_\tau\overset{*}{\rightharpoonup}u\quad\text{in }L^\infty(0,T;H),
\]
\[
\bar u_\tau\to u\quad\text{strongly in }L^2(0,T;H),
\]
\[
\gamma\bar u_\tau\to\gamma u\quad\text{strongly in }L^2(0,T;U),
\]
and
\[
\bar\xi_\tau\rightharpoonup\xi\quad\text{in }L^2(0,T;U^*).
\]
Furthermore,
\[
u'\in L^2(0,T;V^*)\quad(d=2),
\qquad
u'\in L^{4/3}(0,T;V^*)\quad(d=3).
\]
\end{theorem}

\begin{proof}
Proposition~\ref{prop:compact} provides the required compactness, Theorem~\ref{thm:Clarke} identifies the multivalued part, and Lemma~\ref{lem:conv} gives the limit of the convective term.

Let $v\in V$ and $\varphi\in C^3([0,T])$ with $\varphi(T)=0$. We multiply \eqref{eq:BE} by $\tau\varphi_1$, multiply the $n$-th equation \eqref{eq:BDF} by $\tau\varphi_n$, and sum over $n$. The discrete derivative is $S_\tau(v,\varphi)$ and its limit is given by Lemma~\ref{lem:SBP}. The remaining terms pass to the limit by the weak convergence in $L^2(0,T;V)$, the strong convergence of $f$, Lemma~\ref{lem:conv}, and \eqref{eq:xiweak}, respectively. We obtain
\[
\begin{aligned}
&-\int_0^T(u,v)\varphi'\,dt-(u_0,v)\varphi(0)
+\int_0^T\nu(\nabla u,\nabla v)\varphi\,dt\\
&\qquad+\int_0^Tb(u,u,v)\varphi\,dt
+\int_0^T\ip{\xi}{\gamma v}\varphi\,dt
=\int_0^T\ip{f}{v}\varphi\,dt.
\end{aligned}
\]
The identity is in a classical way, first obtained for tensor-product test functions. By density of finite sums of such functions in the corresponding  spaces and the uniform bounds already proven, it extends to all admissible test functions. In particular,
\[
g:=f-Au-B(u,u)-\gamma^*\xi\in L^{q_d}(0,T;V^*)
\]
and the distributional derivative of $u$ is equal to $g$. Hence the evolution equation holds in $V^*$ for a.e. $t\in(0,T)$ and the stated regularity of $u'$ follows.

Moreover, we identify the initial value in $H$. Since
\[
u\in L^\infty(0,T;H),\qquad u'\in L^{q_d}(0,T;V^*)\subset L^1(0,T;V^*),
\]
the standard weak-continuity theorem for the Gelfand triple gives $u\in C_w([0,T];H)$; see \cite[Lemma~III.1.4]{Temam1977}. Lemma~\ref{lem:SBP} shows that $u(0)=u_0$ in $V^*$. Since $H\hookrightarrow V^*$ is injective and both elements belong to $H$, the equality also holds in $H$.
\end{proof}

\begin{corollary}[Energy inequality inherited by the limit]\label{cor:limitenergy}
There exist constants $c_0,C_0,C_1>0$, depending only on the structural data and the Ehrling constants used in the absorption, such that the weak solution obtained in Theorem~\ref{thm:main} satisfies, for every $t\in[0,T]$,
\begin{equation}\label{eq:limitenergy}
\operatorname*{ess\,sup}_{0<s<t}|u(s)|^2
+c_0\int_0^t\|u(s)\|^2\,ds
\le C_0e^{C_1t}
\left(1+|u_0|^2+\int_0^t\|f(s)\|_{V^*}^2\,ds\right).
\end{equation}
In particular, the three-dimensional limit retains the energy-class bound furnished by the discrete scheme.
\end{corollary}

\begin{proof}
We repeat the estimate from Proposition~\ref{prop:energy} up to an arbitrary index $m\le N$. The same discrete Gronwall argument yields
\[
\max_{0\le n\le m}|u^n|^2
+c_0\tau\sum_{n=1}^m\|u^n\|^2
\le C_0e^{C_1t_m}
\left(1+|u^0|^2+\tau\sum_{n=1}^m\|f^n\|_{V^*}^2\right).
\]
For fixed $t\in[0,T]$, choose $m=m(\tau,t)$ such that $t_m\to t$. Jensen's inequality bounds the last discrete sum by $\int_0^{t_m}\|f(s)\|_{V^*}^2ds$. We then pass to the limit using the weak-star lower semicontinuity in $L^\infty(0,t;H)$, the weak lower semicontinuity in $L^2(0,t;V)$, the convergence $u^0_\tau\to u_0$ in $H$, and the weak continuity of $u$. This proves \eqref{eq:limitenergy}.
\end{proof}

\section{Error estimates}\label{sec:error}

We now consider the rate of convergence. This requires additional regularity of the exact solution and the relaxed monotonicity assumption \eqref{eq:relmono},  Error estimates for Rothe approximations of parabolic hemivariational inclusions were obtained in \cite{Kalita2012} and for the incompressible Navier--Stokes equations, second-order BDF estimates were proved in \cite{EmmrichError2004}. We adapt these arguments to the present problem with a Clarke subdifferential term.

The approximation of the right-hand side has to be treated separately. The interval average $f^n$ approximates the value $f(t_n)$ only with first order. Thus the endpoint formulation of \eqref{eq:BDF} with the original loads is, in general, only first-order consistent, although the BDF differentiation formula itself is of second order. We first obtain the corresponding $O(\tau)$ estimate and then modify the righthand side so can improve the result to $O(\tau^2)$.

Throughout this section we assume relaxed monotonicity \eqref{eq:relmono}. We also impose the following regularity hypothesis on one fixed solution of the continuous problem.

\begin{assumption}\label{ass:errorreg}
There exists a solution pair $(u,\xi)$ of \eqref{eq:problem} such that
\begin{equation}\label{eq:errorreg-u}
 u\in C([0,T];V)\cap C^1([0,T];H),
 \qquad
 f\in H^1(0,T;V^*).
\end{equation}
Define the combined parabolic residual
\begin{equation}\label{eq:h-def}
 h:=f-u'.
\end{equation}
For the first-order estimate we assume
\begin{equation}\label{eq:errorreg-h}
 h\in H^1(0,T;V^*).
\end{equation}
We also fix a pointwise selection $\xi:[0,T]\to U^*$ such that
\begin{equation}\label{eq:errorreg-xi}
 \xi(t)\in\partial J(\gamma u(t))\quad\text{for every }t\in[0,T],
\end{equation}
and the evolution equation in \eqref{eq:problem} holds in $V^*$ at every $t\in[0,T]$. No continuity of this selection is assumed. For the second-order estimate we additionally require
\begin{equation}\label{eq:errorreg-h2}
 h\in H^2(0,T;V^*).
\end{equation}
\end{assumption}

\begin{assumption}[Additional regularity for a superconvergent backward-Euler start]\label{ass:startupreg}
When the backward-Euler value $u^1$ from \eqref{eq:BE} is used as the starting value for a second-order estimate, we additionally assume
\begin{equation}\label{eq:startupreg-h}
 h\in W^{1,\infty}(0,T;H).
\end{equation}
\end{assumption}

.

In dimension three the assertions below are conditional on the existence of a solution satisfying Assumption~\ref{ass:errorreg} on the interval under consideration. If such a solution exists only on $[0,T_*]$, all estimates are understood on this shorter interval. 

For a Banach-space-valued function $g$ we write
\[
 (\mathcal A_\tau g)^n:=\frac1\tau\int_{t_{n-1}}^{t_n}g(s)\,ds
\]
and for $n\ge2$ we define
\begin{equation}\label{eq:Qtau}
 (\mathcal Q_\tau g)^n
 :=\frac32(\mathcal A_\tau g)^n-\frac12(\mathcal A_\tau g)^{n-1}.
\end{equation}
A direct calculation gives
\begin{equation}\label{eq:D2Q}
 D_2u(t_n)=(\mathcal Q_\tau u')^n.
\end{equation}
Indeed,
\[
 \frac32\int_{t_{n-1}}^{t_n}u'(s)\,ds
 -\frac12\int_{t_{n-2}}^{t_{n-1}}u'(s)\,ds
 =\frac32u(t_n)-2u(t_{n-1})+\frac12u(t_{n-2}).
\]

\begin{lemma}\label{lem:loadconsistency}
Let $X$ be a Hilbert space.

\begin{enumerate}[label=(\roman*)]
\item If $g'\in L^2(0,T;X)$, then
\begin{equation}\label{eq:firstorder-average}
 \tau\sum_{n=1}^N
 \|(\mathcal A_\tau g)^n-g(t_n)\|_X^2
 \le \frac{\tau^2}{3}\|g'\|_{L^2(0,T;X)}^2.
\end{equation}
\item If $g'\in L^2(0,T;X)$, then
\begin{equation}\label{eq:firstorder-Q}
 \tau\sum_{n=2}^N
 \|(\mathcal Q_\tau g)^n-g(t_n)\|_X^2
 \le C\tau^2\|g'\|_{L^2(0,T;X)}^2.
\end{equation}
\item If $g''\in L^2(0,T;X)$, then
\begin{equation}\label{eq:secondorder-Q}
 \tau\sum_{n=2}^N
 \|(\mathcal Q_\tau g)^n-g(t_n)\|_X^2
 \le C\tau^4\|g''\|_{L^2(0,T;X)}^2,
\end{equation}
where $C$ is independent of $\tau$ and $g$.
\end{enumerate}
\end{lemma}

\begin{proof}
To prove (i), let $a=t_{n-1}$ and $b=t_n$. By Fubini's theorem,
\[
\begin{aligned}
 (\mathcal A_\tau g)^n-g(b)
 &=\frac1\tau\int_a^b\bigl(g(s)-g(b)\bigr)\,ds\\
 &=-\frac1\tau\int_a^b\int_s^b g'(r)\,dr\,ds
 =-\frac1\tau\int_a^b(r-a)g'(r)\,dr.
\end{aligned}
\]
Hence, by the Cauchy--Schwarz inequality,
\[
 \|(\mathcal A_\tau g)^n-g(t_n)\|_X^2
 \le \frac1{\tau^2}\left(\int_a^b(r-a)^2\,dr\right)
      \left(\int_a^b\|g'(r)\|_X^2\,dr\right)
 =\frac{\tau}{3}\int_a^b\|g'(r)\|_X^2\,dr.
\]
After multiplying by $\tau$ and summing with respect to $n$, we obtain \eqref{eq:firstorder-average}.

For (ii), using the definition of $\mathcal Q_\tau$ and part (i) on the two adjacent intervals, we obtain
\[
 \|(\mathcal Q_\tau g)^n-g(t_n)\|_X
 \le \frac32\|(\mathcal A_\tau g)^n-g(t_n)\|_X
 +\frac12\|(\mathcal A_\tau g)^{n-1}-g(t_{n-1})\|_X
 +\frac12\|g(t_{n-1})-g(t_n)\|_X.
\]
Squaring, summing, and using the Cauchy--Schwarz inequality on the last difference gives \eqref{eq:firstorder-Q}.

For (iii), fix $n\ge2$. For $s\in[t_{n-2},t_n]$ we use the representation
\begin{equation}\label{eq:taylor-g}
 g(s)=g(t_n)+(s-t_n)g'(t_n)+R_n(s),
 \qquad
 R_n(s)=\int_{t_n}^{s}(s-r)g''(r)\,dr.
\end{equation}
 The coefficients in $\mathcal Q_\tau$ reproduce constants,
\[
 \frac{3}{2\tau}\int_{t_{n-1}}^{t_n}1\,ds
 -\frac{1}{2\tau}\int_{t_{n-2}}^{t_{n-1}}1\,ds=1,
\]
and $(s-t_n)$ produce zero:
\begin{equation}\label{eq:firstmomentzero}
 \frac{3}{2\tau}\int_{t_{n-1}}^{t_n}(s-t_n)\,ds
 -\frac{1}{2\tau}\int_{t_{n-2}}^{t_{n-1}}(s-t_n)\,ds
 =-\frac{3\tau}{4}+\frac{3\tau}{4}=0.
\end{equation}
Therefore only the remainder $R_n$ contributes. For $s\in[t_{n-2},t_n]$,
\[
 \|R_n(s)\|_X
 \le \left(\int_s^{t_n}|s-r|^2\,dr\right)^{1/2}
      \|g''\|_{L^2(t_{n-2},t_n;X)}
 \le C\tau^{3/2}\|g''\|_{L^2(t_{n-2},t_n;X)}.
\]
Using \eqref{eq:Qtau},
\[
 \|(\mathcal Q_\tau g)^n-g(t_n)\|_X
 \le C\tau^{3/2}\|g''\|_{L^2(t_{n-2},t_n;X)}.
\]
Therefore
\[
 \tau\sum_{n=2}^N\|(\mathcal Q_\tau g)^n-g(t_n)\|_X^2
 \le C\tau^4\sum_{n=2}^N
 \|g''\|_{L^2(t_{n-2},t_n;X)}^2.
\]
Since every point of $(0,T)$ belongs to at most two intervals $(t_{n-2},t_n)$, estimate \eqref{eq:secondorder-Q} follows.
\end{proof}

\begin{lemma}\label{lem:nonlinear-error}
Assume \eqref{eq:relmono}. Let $U,e\in V$, let $\Xi\in\partial J(\gamma U)$ and $\zeta\in\partial J(\gamma(U+e))$. Then
\begin{equation}\label{eq:nonlinear-error-bound}
 \ip{B(U+e,U+e)-B(U,U)}{e}
 +\ip{\zeta-\Xi}{\gamma e}
 \ge -\frac\nu4\|e\|^2-K_d(U)|e|^2,
\end{equation}
where
\[
 K_2(U)\le C(1+\|U\|^2),
 \qquad
 K_3(U)\le C(1+\|U\|^4).
\]
The constants depend on $\nu,m_J,\gamma$ and the Ehrling constants, but not on $U$ or $e$ otherwise.
\end{lemma}

\begin{proof}
Using the skew-symmetry of $b$ and the fact that $U+e$ is divergence free,
\[
\begin{aligned}
 \ip{B(U+e,U+e)-B(U,U)}{e}
 &=b(e,U,e)+b(U,e,e)+b(e,e,e)\\
 &=b(e,U,e).
\end{aligned}
\]
For $d=2$, the Ladyzhenskaya inequality gives us
\[
 |b(e,U,e)|
 \le C\|U\|\,|e|\,\|e\|
 \le \frac\nu8\|e\|^2+C\|U\|^2|e|^2.
\]
For $d=3$, \eqref{eq:3dL4} and Young's inequality with exponents $4/3$ and $4$ give
\[
 |b(e,U,e)|
 \le C\|U\|\,|e|^{1/2}\|e\|^{3/2}
 \le \frac\nu8\|e\|^2+C\|U\|^4|e|^2.
\]
On the other hand, \eqref{eq:relmono} and Lemma~\ref{lem:ehrling}, with the Ehrling parameter chosen sufficiently small, give
\[
 \ip{\zeta-\Xi}{\gamma e}
 \ge -m_J\|\gamma e\|_U^2
 \ge -\frac\nu8\|e\|^2-C|e|^2.
\]
Combining the preceding inequalities gives \eqref{eq:nonlinear-error-bound}.
\end{proof}

We next estimate the first time step,  assumption~\ref{ass:errorreg} yields the accuracy required for the first-order result and for the second-order BDF estimate we need a starting value of order two. Assumption~\ref{ass:startupreg} is a sufficient condition for the backward-Euler step to have this property.

\begin{lemma}\label{lem:startup-error}
Assume \ref{ass:errorreg}, let $u^0=u_0$ exactly, and let $(u^1,\xi^1)$ solve \eqref{eq:BE}.  Then there are $\tau_{\rm st}>0$ and $C>0$, independent of $\tau$, such that for $0<\tau\le\tau_{\rm st}$,
\begin{equation}\label{eq:startup-first-rate}
 |u^1-u(t_1)|^2+\tau\|u^1-u(t_1)\|^2
 \le C\tau^2\|h'\|_{L^2(0,T;V^*)}^2.
\end{equation}
If, in addition, Assumption~\ref{ass:startupreg} holds, then the stronger estimate
\begin{equation}\label{eq:startup-rate}
 |u^1-u(t_1)|\le C\tau^2,
 \qquad
 \tau\|u^1-u(t_1)\|^2\le C\tau^4
\end{equation}
holds.
\end{lemma}

\begin{proof}
Set $U^1=u(t_1)$, $\Xi^1=\xi(t_1)$, and $e^1=u^1-U^1$. Since $u^0=u(0)$, subtracting the continuous equation at $t_1$ from \eqref{eq:BE} gives
\begin{equation}\label{eq:start-error-eq}
 \frac{e^1}{\tau}+Ae^1+
 \bigl(B(u^1,u^1)-B(U^1,U^1)\bigr)
 +\gamma^*(\xi^1-\Xi^1)=\rho^1,
\end{equation}
where
\begin{equation}\label{eq:rho1}
 \rho^1=(\mathcal A_\tau h)^1-h(t_1).
\end{equation}
From Lemma~\ref{lem:loadconsistency}(i), with $X=V^*$, we have
\begin{equation}\label{eq:rho1-vstar}
 \|\rho^1\|_{V^*}^2
 \le \frac\tau3\int_0^\tau\|h'(s)\|_{V^*}^2\,ds.
\end{equation}
We test \eqref{eq:start-error-eq} with $e^1$ and apply Lemma~\ref{lem:nonlinear-error} with $U=U^1$. Thus
\[
 \left(\frac1\tau-K\right)|e^1|^2+\frac{3\nu}{4}\|e^1\|^2
 \le \ip{\rho^1}{e^1},
\]
where $K$ is independent of $\tau$, since $u\in C([0,T];V)$. The right-hand side satisfies
\[
 |\ip{\rho^1}{e^1}|
 \le \|\rho^1\|_{V^*}\|e^1\|
 \le \frac\nu8\|e^1\|^2+C\|\rho^1\|_{V^*}^2.
\]
If $\tau$ is sufficiently small, so that $K\tau\le1/2$, then
\[
 \frac1{2\tau}|e^1|^2+c\|e^1\|^2
 \le C\|\rho^1\|_{V^*}^2.
\]
Combining this inequality with \eqref{eq:rho1-vstar} we obtain
\[
 |e^1|^2+\tau\|e^1\|^2
 \le C\tau^2\int_0^\tau\|h'(s)\|_{V^*}^2\,ds
 \le C\tau^2\|h'\|_{L^2(0,T;V^*)}^2,
\]
and \eqref{eq:startup-first-rate} follows.

Suppose now that Assumption~\ref{ass:startupreg} holds. Then $\rho^1\in H$ and
\begin{align*}
 |\rho^1|
 &=\left|\frac1\tau\int_0^\tau\bigl(h(s)-h(t_1)\bigr)\,ds\right|\\
 &\le \frac1\tau\int_0^\tau\int_s^\tau|h'(r)|\,dr\,ds
 \le \frac\tau2\|h'\|_{L^\infty(0,T;H)}.
\end{align*}
We return to the tested error equation and estimate the residual in the $H$ norm. We get
\[
 \frac1{2\tau}|e^1|^2+c\|e^1\|^2
 \le |\rho^1|\,|e^1|.
\]
It follows that
\[
 |e^1|\le2\tau|\rho^1|\le C\tau^2.
\]
Substitution in the preceding inequality gives $\|e^1\|^2\le C\tau^3$, and consequently
$\tau\|e^1\|^2\le C\tau^4$. This proves \eqref{eq:startup-rate}.
\end{proof}

We shall say that a starting procedure is \emph{second-order accurate} if its value $u^1$ satisfies
\begin{equation}\label{eq:S2-start}
 |u^1-u(t_1)|^2+\tau\|u^1-u(t_1)\|^2\le C_{\rm st}\tau^4
\end{equation}
with a constant independent of $\tau$. By Lemma~\ref{lem:startup-error}, the backward-Euler start has this property under Assumption~\ref{ass:startupreg}.

We now prove the stability estimate which will be used for both choices of the discrete force. We formulate it for an arbitrary sequence $(g^n)$ at the BDF levels.

\begin{theorem}\label{thm:error-stability}
Assume \eqref{eq:relmono} and Assumption \ref{ass:errorreg}.  Let $u^0=u_0$, and for $n\ge2$ let $(u^n,\xi^n)$ satisfy
\begin{equation}\label{eq:generic-bdf-error-scheme}
 D_2u^n+Au^n+B(u^n,u^n)+\gamma^*\xi^n=g^n,
 \qquad \xi^n\in\partial J(\gamma u^n).
\end{equation}
Set
\[
 U^n=u(t_n),\qquad \Xi^n=\xi(t_n),\qquad e^n=u^n-U^n,
\]
and define the consistency residual
\begin{equation}\label{eq:rho-n}
 \rho^n:=g^n-f(t_n)-\bigl(D_2U^n-u'(t_n)\bigr),
 \qquad n\ge2.
\end{equation}
Then there exist $\tau_{\rm err}>0$ and $C>0$, depending on $T$, $\nu$, $m_J$, $\gamma$, the Ehrling constants and $\|u\|_{L^\infty(0,T;V)}$, but independent of $\tau$, such that for $0<\tau\le\tau_{\rm err}$,
\begin{align}
 &\max_{1\le n\le N}|e^n|^2
 +\tau\sum_{n=2}^N\|e^n\|^2
 +\sum_{n=2}^N|e^n-2e^{n-1}+e^{n-2}|^2
 \notag\\
 &\hspace{35mm}\le
 C\left(|e^1|^2+\tau\sum_{n=2}^N\|\rho^n\|_{V^*}^2\right).
\label{eq:generic-error-bound}
\end{align}
The result is valid for both $d=2$ and $d=3$.
\end{theorem}

\begin{proof}
Subtracting the continuous equation at $t_n$ from \eqref{eq:generic-bdf-error-scheme}, we obtain
\begin{equation}\label{eq:error-equation}
 D_2e^n+Ae^n+
 \bigl(B(u^n,u^n)-B(U^n,U^n)\bigr)
 +\gamma^*(\xi^n-\Xi^n)=\rho^n.
\end{equation}
We test the equation with $e^n$. Lemma~\ref{lem:nonlinear-error}, applied with $U=U^n$, $e=e^n$, $\Xi=\Xi^n$, and $\zeta=\xi^n$, produces
\begin{equation}\label{eq:error-nonlinear-combined}
 \ip{B(u^n,u^n)-B(U^n,U^n)}{e^n}
 +\ip{\xi^n-\Xi^n}{\gamma e^n}
 \ge -\frac\nu4\|e^n\|^2-K_u|e^n|^2,
\end{equation}
where one may take
\begin{equation}\label{eq:Ku-dimension}
 K_u\le
 \begin{cases}
 C\bigl(1+\|u\|_{L^\infty(0,T;V)}^2\bigr),&d=2,\\[1mm]
 C\bigl(1+\|u\|_{L^\infty(0,T;V)}^4\bigr),&d=3.
 \end{cases}
\end{equation}
Moreover,
\begin{equation}\label{eq:error-residual-young}
 |\ip{\rho^n}{e^n}|
 \le \|\rho^n\|_{V^*}\|e^n\|
 \le \frac\nu8\|e^n\|^2+C_\nu\|\rho^n\|_{V^*}^2.
\end{equation}
From \eqref{eq:error-equation}--\eqref{eq:error-residual-young} we obtain
\begin{equation}\label{eq:error-before-G}
 (D_2e^n,e^n)+c_\nu\|e^n\|^2
 \le K|e^n|^2+C\|\rho^n\|_{V^*}^2,
\end{equation}
where $c_\nu>0$ and $K>0$ are independent of $\tau$.

We introduce the discrete $G$-energy
\begin{equation}\label{eq:error-G-energy}
 \mathcal F_n:=|e^n|^2+|2e^n-e^{n-1}|^2.
\end{equation}
The identity \eqref{eq:Gidentity} gives
\begin{equation}\label{eq:error-Gidentity}
 4\tau(D_2e^n,e^n)
 =\mathcal F_n-\mathcal F_{n-1}
 +|e^n-2e^{n-1}+e^{n-2}|^2.
\end{equation}
Since $|e^n|^2\le\mathcal F_n$, multiplying \eqref{eq:error-before-G} by $4\tau$ gives
\begin{align}
 \mathcal F_n-\mathcal F_{n-1}
 &+|e^n-2e^{n-1}+e^{n-2}|^2
 +c\tau\|e^n\|^2
 \notag\\
 &\le K_0\tau\mathcal F_n+C_0\tau\|\rho^n\|_{V^*}^2.
\label{eq:error-recursion}
\end{align}
Choose $\tau_{\rm err}>0$ such that $K_0\tau_{\rm err}\le1/2$. Omitting the nonnegative terms and absorbing $K_0\tau\mathcal F_n$, we have
\[
 \mathcal F_n
 \le (1+2K_0\tau)\mathcal F_{n-1}
 +2C_0\tau\|\rho^n\|_{V^*}^2.
\]
The discrete Gronwall lemma gives
\begin{equation}\label{eq:error-E-bound}
 \max_{1\le n\le N}\mathcal F_n
 \le C e^{CT}
 \left(\mathcal F_1+\tau\sum_{n=2}^N\|\rho^n\|_{V^*}^2\right).
\end{equation}
The stability constant cHence \eqref{eq:generic-error-bound} follows.
\end{proof}

We first consider the original scheme with averaged forces.

\begin{theorem}\label{thm:error-first-order}
Assume \ref{ass:errorreg}, \eqref{eq:relmono}, set $u^0=u_0$, and let $(u^n,\xi^n)$ be a solution of the original scheme \eqref{eq:BE}--\eqref{eq:BDF}, with
\[
 f^n=(\mathcal A_\tau f)^n.
\]
Then, for all sufficiently small $\tau$,
\begin{align}
 &\max_{1\le n\le N}|u^n-u(t_n)|
 +\left(\tau\sum_{n=1}^N\|u^n-u(t_n)\|^2\right)^{1/2}
 \notag\\
 &\quad
 +\left(\sum_{n=2}^N
 |(u^n-u(t_n))-2(u^{n-1}-u(t_{n-1}))+(u^{n-2}-u(t_{n-2}))|^2\right)^{1/2}
 \le C\tau,
\label{eq:first-order-rate}
\end{align}
with $C$ independent of $\tau$.
\end{theorem}

\begin{proof}
For $n\ge2$ we apply Theorem~\ref{thm:error-stability} with
\[
 g^n=f^n=(\mathcal A_\tau f)^n.
\]
Using \eqref{eq:D2Q} and $h=f-u'$, we obtain
\begin{align*}
 \rho^n
 &= (\mathcal A_\tau f)^n-f(t_n)
    -\bigl((\mathcal Q_\tau u')^n-u'(t_n)\bigr)\\
 &=\bigl((\mathcal A_\tau f)^n-(\mathcal Q_\tau f)^n\bigr)
   +\bigl((\mathcal Q_\tau h)^n-h(t_n)\bigr).
\end{align*}
The second term is estimated by Lemma~\ref{lem:loadconsistency}(ii). For the first one,
\[
 (\mathcal A_\tau f)^n-(\mathcal Q_\tau f)^n
 =-\frac12\bigl((\mathcal A_\tau f)^n-(\mathcal A_\tau f)^{n-1}\bigr).
\]
Moreover,
\[
 (\mathcal A_\tau f)^n-(\mathcal A_\tau f)^{n-1}
 =\frac1\tau\int_{t_{n-1}}^{t_n}\bigl(f(s)-f(s-\tau)\bigr)\,ds.
\]
By Jensen's and the Cauchy--Schwarz inequalities
\begin{align*}
 \|(\mathcal A_\tau f)^n-(\mathcal A_\tau f)^{n-1}\|_{V^*}^2
 &\le \frac1\tau\int_{t_{n-1}}^{t_n}
       \|f(s)-f(s-\tau)\|_{V^*}^2\,ds\\
 &\le \int_{t_{n-1}}^{t_n}\int_{s-\tau}^{s}
       \|f'(r)\|_{V^*}^2\,dr\,ds\\
 &\le \tau\int_{t_{n-2}}^{t_n}\|f'(r)\|_{V^*}^2\,dr.
\end{align*}
Since the intervals $(t_{n-2},t_n)$ have overlap at most twice, Lemma~\ref{lem:loadconsistency}(ii) applied to $h$ gives
\begin{equation}\label{eq:first-residual}
 \tau\sum_{n=2}^N\|\rho^n\|_{V^*}^2
 \le C\tau^2\left(\|f'\|_{L^2(0,T;V^*)}^2
 +\|h'\|_{L^2(0,T;V^*)}^2\right).
\end{equation}
By Lemma~\ref{lem:startup-error},
$|e^1|^2+\tau\|e^1\|^2\le C\tau^2$ under Assumption~\ref{ass:errorreg}. Inserting this and \eqref{eq:first-residual} into \eqref{eq:generic-error-bound}, we obtain
\[
 \max_n|e^n|^2
 +\tau\sum_{n=1}^N\|e^n\|^2
 +\sum_{n=2}^N|e^n-2e^{n-1}+e^{n-2}|^2
 \le C\tau^2.
\]
Hence we have proven \eqref{eq:first-order-rate}.
\end{proof}

\subsection{Higher order convergence}
We define
\begin{equation}\label{eq:fhat}
 \widehat f^1:=f^1,
 \qquad
 \widehat f^n:=\frac32f^n-\frac12f^{n-1}
 =(\mathcal Q_\tau f)^n,
 \qquad n\ge2.
\end{equation}
We replace $f^n$ in \eqref{eq:BDF} by $\widehat f^n$ for $n
\ge 2$, cf. \cite{EmmrichError2004}. The solvability argument is unchanged. Moreover,
\[
 \tau\sum_{n=2}^N\|\widehat f^n\|_{V^*}^2
 \le C\tau\sum_{n=1}^N\|f^n\|_{V^*}^2
 \le C\|f\|_{L^2(0,T;V^*)}^2,
\]
and therefore the a priori estimates from Section~3 remain valid, possibly with different constants.

\begin{theorem}\label{thm:error-second-order}
Assume \ref{ass:errorreg}, \eqref{eq:errorreg-h2}, and \eqref{eq:relmono}.  Let $u^0=u_0$ and suppose the starting value $u^1$ is second-order accurate in the sense of \eqref{eq:S2-start}.  For $n\ge2$ let $(u^n,\xi^n)$ solve
\begin{equation}\label{eq:BDF-secondorderload}
 D_2u^n+Au^n+B(u^n,u^n)+\gamma^*\xi^n=\widehat f^n,
 \qquad
 \xi^n\in\partial J(\gamma u^n),
\end{equation}
with \eqref{eq:fhat}.  Then, for all sufficiently small $\tau$,
\begin{align}
 &\max_{1\le n\le N}|u^n-u(t_n)|
 +\left(\tau\sum_{n=1}^N\|u^n-u(t_n)\|^2\right)^{1/2}
 \notag\\
 &\quad
 +\left(\sum_{n=2}^N
 |(u^n-u(t_n))-2(u^{n-1}-u(t_{n-1}))+(u^{n-2}-u(t_{n-2}))|^2\right)^{1/2}
 \le C\tau^2.
\label{eq:second-order-rate}
\end{align}
The constant $C$ is independent of $\tau$.  In particular, the backward-Euler start \eqref{eq:BE} is admissible when Assumption \ref{ass:startupreg} holds, by Lemma~\ref{lem:startup-error}.  In dimension three the assertion is conditional on the existence of the regular solution from Assumption \ref{ass:errorreg} on the time interval considered; if such a solution is known only on $[0,T_*]$, the estimate holds on that interval.
\end{theorem}

\begin{proof}
For $n\ge2$ we apply Theorem~\ref{thm:error-stability} with
\[
 g^n=\widehat f^n=(\mathcal Q_\tau f)^n.
\]
By \eqref{eq:D2Q} and the definition of $h$, the two consistency errors combine as follows:
\begin{align*}
 \rho^n
 &=\bigl((\mathcal Q_\tau f)^n-f(t_n)\bigr)
   -\bigl((\mathcal Q_\tau u')^n-u'(t_n)\bigr)\\
 &=(\mathcal Q_\tau h)^n-h(t_n).
\end{align*}
 Lemma~\ref{lem:loadconsistency}(iii) gives
\begin{equation}\label{eq:second-residual}
 \tau\sum_{n=2}^N\|\rho^n\|_{V^*}^2
 \le C\tau^4\|h''\|_{L^2(0,T;V^*)}^2.
\end{equation}
Using the estimate \eqref{eq:S2-start} together with \eqref{eq:second-residual} in \eqref{eq:generic-error-bound}, we get
\[
 \max_n|e^n|^2
 +\tau\sum_{n=1}^N\|e^n\|^2
 +\sum_{n=2}^N|e^n-2e^{n-1}+e^{n-2}|^2
 \le C\tau^4.
\]
This proves \eqref{eq:second-order-rate}.
\end{proof}


\begin{thebibliography}{99}

\bibitem{FangHanMigorskiSofonea2016}
C. Fang, W. Han, S. Mig\'orski and M. Sofonea,
\emph{A class of hemivariational inequalities for nonstationary Navier--Stokes equations},
Nonlinear Analysis: Real World Applications 31 (2016), 257--276.


\bibitem{BartoszSzafraniecZhao2021}
K. Bartosz, P. Szafraniec and J. Zhao,
\emph{Convergence of a double step scheme for a class of parabolic Clarke subdifferential inclusions},
Communications in Nonlinear Science and Numerical Simulation 102 (2021), 105940.


\bibitem{BartoszSzafraniec2024}
K. Bartosz and P. Szafraniec,
\emph{Convergence of a double step scheme for a class of second order Clarke subdifferential inclusions},
Nonlinear Analysis: Real World Applications 78 (2024), 104092.


\bibitem{Emmrich2004}
E. Emmrich,
\emph{Stability and convergence of the two-step BDF for the incompressible Navier--Stokes problem},
International Journal of Nonlinear Sciences and Numerical Simulation 5 (2004), 199--209.

\bibitem{Dahlquist1978}
G. Dahlquist,
\emph{$G$-stability is equivalent to $A$-stability},
BIT 18 (1978), 384--401.

\bibitem{NevanlinnaOdeh1981}
O. Nevanlinna and F. Odeh,
\emph{Multiplier techniques for linear multistep methods},
Numerical Functional Analysis and Optimization 3 (1981), 377--423.


\bibitem{MigorskiOchalSofonea2013}
S. Mig\'orski, A. Ochal and M. Sofonea,
\emph{Nonlinear Inclusions and Hemivariational Inequalities: Models and Analysis of Contact Problems},
Advances in Mechanics and Mathematics, vol. 26, Springer, New York, 2013.

\bibitem{EmmrichError2004}
E. Emmrich,
\emph{Error of the two-step BDF for the incompressible Navier--Stokes problem},
ESAIM: Mathematical Modelling and Numerical Analysis 38 (2004), 757--764.


\bibitem{Kalita2012}
P. Kalita,
\emph{Regularity and Rothe method error estimates for parabolic hemivariational inequality},
Journal of Mathematical Analysis and Applications 389 (2012), 618--631.


\bibitem{Simon1987}
J. Simon,
\emph{Compact sets in the space $L^p(0,T;B)$},
Annali di Matematica Pura ed Applicata 146 (1987), 65--96.



\bibitem{AubinCellina1984}
J.-P. Aubin and A. Cellina,
\emph{Differential Inclusions: Set-Valued Maps and Viability Theory},
Springer-Verlag, Berlin, 1984.

\bibitem{Clarke1990}
F. H. Clarke,
\emph{Optimization and Nonsmooth Analysis},
SIAM, Philadelphia, 1990.


\bibitem{Temam1977}
R. Temam,
\emph{Navier--Stokes Equations: Theory and Numerical Analysis},
North-Holland, Amsterdam, 1977.

\end{thebibliography}
\end{document}